\documentclass[oneside, reqno]{amsart}
\usepackage{amsmath,amssymb,latexsym,amsfonts,enumerate,longtable,amscd}
\usepackage{setspace}
\usepackage{xy}
\usepackage{paralist,xy}

\everymath{\displaystyle} 
\input xy
\xyoption{all}

\newtheorem{thm}{Theorem}[section]
\newtheorem{corollary}[thm]{Corollary}
\newtheorem{lemma}[thm]{Lemma}
\newtheorem{question}[thm]{Question}
\newtheorem{definition}[thm]{Definition}
\newtheorem{remark}[thm]{Remark}
\newtheorem{example}[thm]{Example}

\newtheorem{proposition}[thm]{Proposition}
\numberwithin{equation}{section}

\numberwithin{equation}{subsection}
\begin {document}
\title{CLASSES OF SOME HYPERSURFACES IN THE GROTHENDIECK RING OF
VARIETIES}
\author{Emel Bilgin}
\maketitle
\begin{abstract}
Let $X$ be a projective hypersurface in $\mathbb{P}_k^n$ of degree
$d \leq n$. In this paper we study the relation between the class
$[X]$ in $K_0(Var_k)$ and the existence of $k$-rational points.
Using elementary geometric methods we show, for some particular $X$,
that $X(k)\neq \emptyset $ if and only if $[X] \equiv 1$ modulo
$\mathbb{L}$ in $K_0(Var_k)$. More precisely we consider the
following cases: a union of hyperplanes, a quadric, a cubic
hypersurface with a singular $k$-rational point, and a quartic which
is a union of two quadrics one of which being smooth.
\end{abstract}
\tableofcontents
\section{Introduction}
Let $k$ be a field, let $Var_k$ denote the category of varieties
over $k$. The Grothendieck ring of varieties, denoted by
$K_0(Var_k)$, is a ring whose first appearance in the literature
dates back to 1964, in a letter of A. Grothendieck to J. P. Serre
\cite{SeCo}. Even though this ring has been studied from many
aspects since then, still only little is known about it. In the
article \cite{Poonen}, B. Poonen proves that if the base field $k$
is of characteristic zero, then $K_0(Var_k)$ is not a domain, by
providing some example of zero divisors. Another example for the
zero divisors of $K_0(Var_k)$ was given by J. Koll\'{a}r in
\cite{Ko}. F. Heinloth \cite{F.H} defined $K_0(Var_k)$ in a
different way, introducing blow up relations on smooth projective
varieties, and proved that this definition is equivalent to the
classical one. M. Larsen and V. A. Lunts \cite{LaLu} established an
interesting relation of $K_0(Var_k)$ with stable rationality, which
is a notion weaker than rationality.

    Let $X$ be a hypersurface of degree $d$ in $\mathbb{P}_k^n$, where
$d\leq n$. Let $[X]$ denote its class in $K_0(Var_k)$. The main
question that concerns us is whether it is true that $X(k)$ is
nonempty if and only if
\begin{equation}\label{intro_main_equiv}[X]\equiv 1 \hspace{.3cm}\text{ modulo
}\hspace{.3cm}[\mathbb{A}_k^1]\end{equation} in $K_0(Var_k)$, or
not. We explain and prove that for some certain types of
hypersurfaces this question has a positive answer.

    After the completion of this thesis, we came to know of an
independent work submitted on Arxiv by X. Liao \cite{kil}, which has
a significant intersection with what we have done. Whereas our main
question is based on the relation between the existence of a
$k$-rational point and the class in $K_0(Var_k)$, with $k$ an
arbitrary field, in \cite{kil} $k$ is assumed to be algebraically
closed of characteristic zero and the focus is on the
"$\mathbb{L}$-rationality". Other than this, the main parts of this
work that differ from \cite{kil} are the following theorems where we
prove that Question \ref{question} is positively answered: Theorem
\ref{hyp2}, which is about a hypersurface $X$ over a field $k$ such
that $X\times_k L$ is a union of $d\leq n$ hyperplanes where $L / k$
is a finite Galois extension, and Theorem \ref{2quadrics}, where we
study the union of two quadrics, one of which is assumed to be
smooth, over an algebraically closed field $k$ of characteristic
zero. \subsection*{Acknowledgements}
    I would like to thank my Ph.D. thesis advisor Prof. Dr. H\'{e}l\`{e}ne Esnault
for giving me the subject of this thesis and for her great support
and encouragement, and also for creating the chance of being in such
a rich mathematical environment. I am greatly indebted to Dr. Andre
Chatzistamatiou for his very valuable guidance and for the many
instructive discussions: most of the ideas consisting this work come
partly from him. I am deeply grateful also for all the time and
energy he spent for the careful reading of my work. I would like to
thank Prof. Dr. Georg Hein for the useful hints, in particular for
the construction of Example \ref{nonrat_sing}. Finally I would like
to thank Sonderforschungsbereich / Transregio 45 for the financial
support during the preparation of this work.
\section{Grothendieck Ring of Varieties} \label{halka-denge}  Let us recall some
definitions and properties that will be used.
\begin{definition} {\rm \cite[Definition 2.1]{Nic}} \label{K_0(Var)}
Let $k$ be a field. Let $Var_k$ denote the category of varieties
over $k$. Note that what we mean here by a variety over $k$ is a
reduced separated scheme of finite type over $k$. The Grothendieck
group $K_0(Var_k)$ is the abelian group generated by isomorphism
classes of varieties over $k$, with the relation $$
\hspace{1.1cm}[X]= [X \setminus Y] + [Y]\text{,}\hspace{1cm} $$ if
$Y \subset X $ is a closed subvariety. The  product
$$[X]\cdot[Z]:=[X\times_k Z]$$ defines a ring structure on
$K_0(Var_k)$.
\end{definition}
Note that the zero of this ring is $0:= [\emptyset]$, and the
multiplicative identity is $1:=[Spec(k)]$. We denote by $\mathbb{L}$
the class of $\mathbb{A}_k^1$ in $K_0(Var_k)$.

\begin{remark}{\rm
We could actually drop the requirement of reducedness in Definition
\ref{K_0(Var)}. Let $X$ be a separated scheme of finite type over
$k$, and let $X_{red}$ denote the reduced scheme associated to $X$.
Applying the scissor relations given in Definition \ref{K_0(Var)} to
the natural closed immersion $X_{red}\rightarrow X$ we obtain
$$[X]=[X\setminus X_{red}]+[X_{red}]=[X_{red}]\text{.}$$
}
\end{remark}
Let $\pi$ : $X\rightarrow Y$ be a morphism of varieties over $k$.
Recall that $\pi$ is said to be a Zariski locally trivial fibration
with fiber $F$ if each closed point $y\in Y$ has a Zariski open
neighborhood $U$ such that the pre-image $\pi^{-1}(U)$ is isomorphic
over $k$ to $F\times_k U$. We will frequently use the following
remark for our calculations in $K_0(Var_k)$.
\begin{remark} \label{trivial_fibration} {\rm
A Zariski locally trivial fibration $\pi : X\rightarrow Y$ becomes
trivial in the Grothendieck ring of varieties, i.e. one gets
$[X]=[F]\cdot[Y]$. This follows from the defining relations of
$K_0(Var_k)$, using induction on the number of open neighborhoods
$U_i$ in $Y$ which has $\pi^{-1}(U_i)=F\times_k U_i$, that covers
$Y$. }
\end{remark}
 Consider the blow up of a smooth projective
variety $X$ along a closed smooth subvariety $Z \subset X$ of
codimension $r$. For such a blow up we also have $$[Bl_Z X]
-[E]=[X]-[Z]$$ in $K_0(Var_k)$. This equation actually yields
another definition for the Grothendieck ring of varieties, which is
proven to be an equivalent presentation of $K_0(Var_k)$, by F.
Heinloth in her paper \cite{F.H}.
\begin{remark} {\rm (\cite[Theorem 3.1]{F.H})\label{Bittner_1} Let
$k$ be a field of characteristic zero, and let $K_0^{bl}(Var_k)$
denote the abelian group generated by the isomorphism classes of
smooth projective varieties with the relations $[\emptyset]_{bl}=0$
and $[Bl_Z X]_{bl} -[E]_{bl}=[X]_{bl}-[Z]_{bl}$ whenever $Bl_Z X$ is
the blow up of a smooth projective variety $X$ along a smooth closed
subvariety $Z$ with the exceptional divisor $E$. Then F. Heinloth
proves that the ring homomorphism
\begin{eqnarray*}
K_0^{bl}(Var_k) & \longrightarrow & K_0(Var_k) \\
{[X]_{bl}} &\longmapsto & {[X]}
\end{eqnarray*} is an isomorphism.}
\end{remark}
\subsection{Stable Birationality} \label{stable_rationality}
\begin{definition}
Let $X, Y $ be irreducible varieties over a field $k$. Then $X$ and
$Y$ are called stably birational if there exists integers $m,n \geq
0$ such that $X\times_k \mathbb{P}_k^{n}$ is birational to
$Y\times_k \mathbb{P}_k^m$ over $k$.
\end{definition}

     Note that stable birationality defines an equivalence
relation $\sim_{SB}$ on the set of smooth complete irreducible
varieties over $k$. Let $SB$ denote the set of equivalence classes
of this relation. We will denote by $\mathbb{Z}[SB]$ the free
abelian group generated by $SB$.

    In the article of M. Larsen and V. A. Lunts \cite{LaLu}, for a
smooth complete irreducible variety it is proven that stable
rationality is a necessary and sufficient condition for being
equivalent to $1$ modulo the class of affine line in
$K_0(Var_{\mathbb{C}})$. Let us recall this important result.
\begin{thm}\label{larsen_lunts_1}{\rm{\cite[Theorem 2.3]{LaLu}} }
Let $G$ be an abelian commutative monoid and $\mathbb{Z}[G]$ be the
corresponding monoid ring. Denote by $\mathcal{M}$ the
multiplicative monoid of isomorphism classes of smooth complete
irreducible varieties over $\mathbb{C}$. Let
$$\Psi \text{ : } \mathcal{M} \longrightarrow G$$ be a homomorphism
of monoids such that
\begin{enumerate}
 \item $\Psi([X])= \Psi([Y])$ if $X$ and $Y$ are birational;
  \item $\Psi([\mathbb{P}_{\mathbb{C}}]) = 1$  for all $n \geq 0$.
\end{enumerate}
Then there exists a unique ring homomorphism $$\Phi \text{ : } K_0(Var_{\mathbb{C}}) \longrightarrow \mathbb{Z}[G]$$ such that
$\Phi([X])=\Psi([X])$ for $[X]\in \mathcal{M}$.
\end{thm}
Now taking $G=SB$ in Theorem \ref{larsen_lunts_1} one gets the ring
homomorphism
 $\Phi_{SB}$ : $K_0(Var_{\mathbb{C}})\rightarrow \mathbb{Z}[SB]$,
 which is clearly surjective.
\begin{proposition} {\rm{\cite[Proposition 2.7]{LaLu}}}
\label{larsen_lunts_2}With the assumptions and notations of Theorem
\ref{larsen_lunts_1}, the kernel of the ring homomorphism
$\Phi_{SB}$ defined above is the principal ideal generated by the
class $\mathbb{L}$.
\end{proposition}
\begin{remark}{\rm \label{remark_ll2}{
For $X$ connected smooth projective, Proposition
\ref{larsen_lunts_2} implies that $X$ is stably rational if and only
if $[X] \equiv 1 \text{ mod } \mathbb{L}$ in
$K_0(Var_{\mathbb{C}})$. Note that however we have no longer this
implication for $X$ singular or non-complete in general. This is
because $\Phi_{SB}$, by construction, sends only the class of a
smooth complete variety to its own equivalence class in
$\mathbb{Z}[SB]$, i.e. for $X$ not smooth complete, the equivalence
$X\sim_{SB} Spec(\mathbb{C})$ does not necessarily imply that
$\Phi_{SB}([X])=1$.} }
\end{remark}

    As it is pointed out by J. Koll\'{a}r in \cite{Ko}, the result of
Larsen and Lunts holds in fact over any field of characteristic
zero: this follows from the presentation of $K_0(Var_k)$ by means of
the blow up relations (see Remark \ref{Bittner_1}).

    The following proposition, which follows directly from Theorem
\ref{larsen_lunts_1}, demonstrates a noteworthy connection between
$K_0(Var_k)$ and $k$-rational points.
\begin{proposition}\label{proposition_ll3}
Let $k$ be a field of characteristic zero, let $X$ be a smooth
connected complete variety of dimension $m$ over $k$. If $[X]\equiv
1 \text{ }\hspace{.1cm} \text{mod } \text{ }\mathbb{L}$, then
Proposition \ref{larsen_lunts_2} implies that $X(k)\neq \emptyset$.
\end{proposition}
\begin{proof}
As it is mentioned in Remark \ref{remark_ll2}, we know that $X$ is
stably rational, i.e. there is an $a \geq 0$ such that $X \times
\mathbb{P}_k^a$ is rational. Let $\phi :
\mathbb{P}_k^{m+a}\dashrightarrow X \times \mathbb{P}_k^a$ be a
birational map, and let $U\subset \mathbb{P}_k^{m+a}$ be the open
subset on which the map $\phi$ is defined, and is an isomorphism. It
is then guaranteed that $U(k)\neq \emptyset $ since we assume
$char(k)=0$, which implies that $k$ is an infinite field. Hence $X
\times \mathbb{P}_k^a$ has a $k$-rational point, which implies that
$X$ has a $k$-rational point.
\end{proof}
    As an effort to make the point of Remark \ref{remark_ll2} clearer,
let us illustrate some negative examples.  For non-complete (stably)
rational varieties, consider simply the affine space
$\mathbb{A}_k^n$. Although it is rational, one has the equivalence
$[\mathbb{A}_k^n] \equiv 0$ mod $\mathbb{L}$ in $K_0(Var_k)$, since
$[\mathbb{A}_k^n]= \mathbb{L}^n$, for all $n\geq 0$.

    Below we give an easy example of singular varieties whose
equivalence classes are $1$ mod $\mathbb{L}$ even though they are
not stably rational.
\begin{example}\label{cone} Let $X \subset \mathbb{P}_k^n$ be a nonrational cone
over a hypersurface $Z\subset \mathbb{P}_k^{n-1}$. Note that a cone
is always singular, and it always has a $k$-rational point.
Therefore after a change of coordinates, we can assume that
$x:=[0:\dots:0:1]\in X$. The projection
\begin{eqnarray*} \mathbb{P}_k^n \setminus \{x\}
&\longrightarrow&
\mathbb{P}_k^{n-1} \\
{[x_0:\dots:x_n]} &\longmapsto& {[x_0:\dots:x_{n-1}]}
\end{eqnarray*}
is a Zariski locally trivial affine fibration. The following is a
cartesian diagram:
$$
\xymatrix{
 \mathbb{P}_k^n \setminus \{x\}\ar[r] &\mathbb{P}_k^{n-1} \\
X \setminus \{x\}\ar[r] \ar@{^{(}->}[u]& Z  \ar@{^{(}->}[u]   }
$$
Thus we get that $X \setminus \{x\} \rightarrow Z$ is also a Zariski
locally trivial affine fibration, which implies by Remark
\ref{trivial_fibration} that $[X \setminus \{x\}]= \mathbb{L} \cdot
[Z]$. Hence
$$[X]=1+\mathbb{L} \cdot [Z] \text{,}$$ and we see that $[X]\equiv
1$ modulo $\mathbb{L}$ in $K_0(Var_k)$.
\end{example}
    Now let us consider the case of a stably rational variety that is
not smooth. In the below examples one sees that such a variety need
not have its class in $K_0(Var_k)$ equivalent to $1$ modulo
$\mathbb{L}$.
\begin{example}
Let $X \subset \mathbb{P}_{\mathbb{C}}^2$ be the projective nodal
curve, more precisely the projective hypersurface given by the
homogenization of the equation $$y^2=x^3+x^2\text{.}$$ It is well
known that $X$ is a rational curve with a unique singular point at
the origin $O$. Let $\pi:\widetilde{X}\rightarrow X$ be the blow up
of $X$ at $O$, and let $E:=\pi^{-1}(O)$. Since the exceptional
divisor $E$ consists of two disjoint $\mathbb{C}$-rational points,
we have $[E]=2$ in $K_0(Var_k)$. On the other hand, the blow up
$\widetilde{X}$ is the normalization of the curve $X$, and therefore
it is a projective smooth rational irreducible variety. Hence Remark
\ref{remark_ll2} implies that $[\widetilde{X}]\equiv 1$ mod
$\mathbb{L}$ in $K_0(Var_k)$. Thus we obtain
\begin{eqnarray*}[X] &=& [X\setminus O]+[O] \\ &=&[\widetilde{X}\setminus E]+[O]
\\ &=&[\widetilde{X}] - [E] + [O] \\ &\equiv & 0 \text{
mod } \mathbb{L}\text{.}\end{eqnarray*}
\end{example}
    The following example illustrates that even a normal singular
variety, i.e. a variety with singularities at least codimension $2$,
might have its class in $K_0(Var_k)$ not equivalent to $1$ modulo
$\mathbb{L}$.
\begin{example} \label{nonrat_sing}
Let $X$ be a rational normal projective surface over an
algebraically closed field $k$ with a unique singular point $x \in
X$, and let $\widetilde{X} \rightarrow X $ be the minimal resolution
of singularities of $X$ where $\widetilde{X}$ is the blow up of $X$
at $x$. Let moreover the exceptional divisor $E:=\pi^{-1}(x)$ be a
union of three lines in $\mathbb{P}_k^2$ which form a cycle, i.e.
the intersection of all three lines is empty. Note that one has
indeed such a surface, and below we will sketch its construction.
But let us first show that $[X] \not\equiv 1$ mod $\mathbb{L}$ in
$K_0(Var_k)$. Being the blow up of a projective rational variety,
$\widetilde{X}$ is also rational and projective, since it is smooth
we can conclude by Remark \ref{remark_ll2} that $[\widetilde{X}]
\equiv 1$ mod $\mathbb{L}$ in $K_0(Var_k)$. Also note that
$$[E]=3[\mathbb{P}_k^1]- 3= 3 \mathbb{L}\equiv 0 \text{ mod } \mathbb{L}$$
in $K_0(Var_k)$. Thus we get
\begin{eqnarray*}
[X] &=& [\widetilde{X} \setminus E]+ [\{x\}] \\
    &=& [\widetilde{X}] - [E] + 1 \\
    &\equiv& 2-[E] \text{ mod } \mathbb{L} \\
    &\equiv& 2 \text{ mod } \mathbb{L}
\end{eqnarray*}
in $K_0(Var_k)$.

    {\rm Let us now roughly describe the construction of $X$ of Example
\ref{nonrat_sing}. Let $L_1,L_2,L_3\subset \mathbb{P}_k^2$ be
distinct lines such that the intersection $L_1\cap L_2\cap L_3$ is
empty, so that they form a cycle. Now we mark four distinct points
$\{p_{i1},p_{i2},p_{i3},p_{i4}\}=:P_i$ on each $L_i$, where none of
the $p_{ij}$'s is the intersection point of $L_i$ and $L_j$. Write
$$S:= P_1 \cup P_2 \cup P_3\text{.}$$ Let $\widetilde{X}:= Bl_S
\mathbb{P}_k^2 \rightarrow \mathbb{P}_k^2$ be the blow up of
$\mathbb{P}_k^2$ at $S$.  We claim that $$ R:=L_1\cup L_2 \cup L_3
\subset \widetilde{X} $$ is contractible, i.e. there exists a
birational morphism
$$f: \widetilde{X} \longrightarrow X $$ such that
$f(R)=\{x\}$, where $x$ is a point of $X$, and $f: \widetilde{X}
\setminus R\rightarrow X\setminus \{x\}$ is an isomorphism. Recall
that to be contractible a curve needs to be negative definite which
means that the matrix of the irreducible components of the curve is
negative definite: in our example the matrix is \vspace{.2cm}
\begin{equation*}\left(L_i\cdot L_j\right)_{i,j} = \left(\begin{matrix}
-3 &  1  & 1 \\
1 & -3  & 1\\
1 & 1 & -3  \\
\end{matrix}\right)\vspace{.3cm}\end{equation*} since $L_i\cdot L_i =
-3$ for all $i$, and $L_i\cdot L_j=1$ for all $i\neq j$. Thus $R$ is
negative definite.

    Denote by $E_{ij}$ the exceptional divisor of $p_{ij}$ in
$\widetilde{X}$. Let $H$ be a general line in $\mathbb{P}_k^2$.
Consider the divisor $D=4H - \displaystyle \sum_{i=1}^3
\displaystyle \sum_{j=1}^4 E_{ij}$ on $\widetilde{X}$. One can prove
that the complete linear system $|mD|$ is base point free for a
sufficiently large integer $m
>> 0$, and therefore $|mD|$ defines a morphism $$f: \widetilde{X} \longrightarrow
\mathbb{P}_k^N \text{.}$$ Note that we have $D\cdot L_i=4H\cdot
L_i-\displaystyle \sum_{j=1}^4 E_{ij}\cdot L_i=0$ for all $i$.
Besides $D^2=16-12=4>0$, and $D \cdot E_{ij}= -E_{ij} \cdot E_{ij} =
1$. Let $C \subset Bl_S \mathbb{P}_k^2$ be any irreducible curve
with $C \not \in \{L_i, E_{ij}\}$ for all $i,j$. Then we can
identify it with the irreducible curve it comes from in the
projective plane $\mathbb{P}_k^2$. Hence one gets $$D\cdot C = 4
deg(C)- \displaystyle \sum_{i=1}^3 \displaystyle \sum_{j=1}^4
mult_{p_{ij}}(C)$$ where $mult_{p_{ij}}(C)$ denotes the multiplicity
of the curve $C$ at the point $p_{ij}$. By the definition of
multiplicity one has $mult_{p_{ij}}(C \cap L_i) \geq
mult_{p_{ij}}(C)$, which yields
\begin{eqnarray*}
D\cdot C &=& 4deg(C)-\sum_{i=1}^3\sum_{j=1}^4 mult_{p_{ij}}(C) \\
&\geq&
4deg(C)- \sum_{i=1}^3\sum_{j=1}^4 mult_{p_{ij}}(C\cap L_i) \\
&\geq &4deg(C)- 3deg(C) \\
&=& {deg(C)} >0 \text{.}
\end{eqnarray*}
    Thus $f$ contracts precisely $R$, and the image $X:=
f(\widetilde{X})$ has the desired properties of Example
\ref{nonrat_sing}.}
\end{example}
\section{Introduction to The General Conjecture} \label{soru-motif}  Let $k$
be a field of characteristic zero, let $X$ be a hypersurface over
$k$ of degree $d$ in $\mathbb{P}_k^n$, where $d \leq n$.
\begin{remark}{\rm Let us recall that a smooth hypersurface $X \subset
\mathbb{P}_k^n$ of degree $d$ is a Fano variety, i.e. its canonical
bundle is anti-ample if and only if $d\leq n$.} \end{remark}
\subsection{Category of Motives} \cite[Chapter 4]{Motif} \noindent Let $k$ be a field, and
let $\mathcal{V}_k$ denote the category of smooth projective schemes
over $k$. For $X \in \mathcal{V}_k$, and for a non-negative integer
$d \leq dimX$, recall that the group of codimension $d$ cycles
$Z^d(X)$ is the free abelian group generated by the codimension $d$
subvarieties of $X$. Then the Chow group of $X$ of codimension $d$
is defined as $CH^d(X): = Z^d(X) / \sim $ , where $\sim$ denotes the
rational equivalence. Note that we consider $CH^d(X)$ to be with
$\mathbb{Q}$-coefficients, unless otherwise stated.
\begin{definition}
Let $X,Y\in \mathcal{V}_k$, and let $X_i $ be the connected
components of $X$. The group of correspondences of degree $r$ from
$X$ to $Y$ is defined to be $$Hom^r(X,Y):= \bigoplus_i CH^{dim
X_i+r}(X_i\times Y)\text{.}$$ We will simply denote by $Hom(X,Y)$
the group of the correspondences of degree zero.
\end{definition}
\noindent Let $P\in Hom(X,Y)$, $Q\in Hom(Y,Z)$. Then the composition
is given by
$$Q\circ P:= (p_{13})_*(p_{12}^*(P)\cdot p_{23}^*(Q))$$ where $p_{ij}$
are the projection maps from $X\times Y\times Z$ to the product of
the $i$-th and $j$-th factors, and $(p_{ij})_*$ , $p_{ij}^*$ are the
pull-back and push-forward of $p_{ij}$ in the Chow groups,
respectively. Note also that $\cdot$ denotes the intersection
product in $CH^{dim X}(X\times Y\times Z)$. This composition gives
$Hom(X,X)$ a $\mathbb{Q}$-algebra structure.
\begin{definition}
A pair $(X,P)$ with $X \in \mathcal{V}_k$ and $P\in Hom(X,X)$ a
projector, i.e. $P=P\circ P$, is called a motive.
\end{definition}
\noindent Motives form a category denoted by $\mathcal{M}_k$ with
morphism groups
$$Hom_{\mathcal{M}_k}((X,P),(Y,Q)):= Q \circ Hom(X,Y)\circ P \subset Hom(X,Y) \text{.}$$
Note that the identity morphism of a motive $(X,P)$ is the projector
$P$. The sum and the product in $\mathcal{M}_k$ are defined by
disjoint union and product:
\begin{eqnarray*}
(X,P)\oplus (Y,Q) &=& (X\sqcup Y, P+Q) \\
(X,P)\otimes (Y,Q)&=& (X\times Y, P\times Q) \text{.}
\end{eqnarray*}
The motive associated with $X \in \mathcal{V}_k$ is the motive $(X,
id_X)$ where $id_X \in Hom(X,X)$ is the class of the diagonal
$\Delta_X$ in the Chow group. There is a functor $$h\text{ : }
\mathcal{V}_k^{op} \longrightarrow \mathcal{M}_k$$ given on objects
by $h(X):=(X,id_X)$, and on morphisms by $h(\varphi):=
[^t\Gamma_{\varphi}]$, the class of the transpose of the graph of
$\varphi: Y\rightarrow X$. Note that $\mathbb{Q}(0):=h(Spec(k))$ is
the identity for the product and it is called the unit motive.
Consider the motive $h(\mathbb{P}^1_k)$. For any $x\in
\mathbb{P}_k^1(k)$ we have
\begin{equation}\label{P1decomp}[\Delta_{\mathbb{P}_k^1}]=[\{x\}\times
{\mathbb{P}_k^1}]+[{\mathbb{P}_k^1} \times \{x\}]\end{equation} in
$CH^{1}(\mathbb{P}_k^1\times \mathbb{P}_k^1)$. Hence
\begin{equation}\label{lefschetz_decomp}
h({\mathbb{P}_k^1})=({\mathbb{P}_k^1},{\{x\}\times\mathbb{P}_k^1})\oplus({\mathbb{P}_k^1},{\mathbb{P}_k^1}\times
\{x\})\text{.}\end{equation} The first summand of the decomposition
\ref{lefschetz_decomp} is isomorphic to $\mathbb{Q}(0)$, and the
latter summand $$({\mathbb{P}_k^1}, {\mathbb{P}_k^1}\times
\{x\})=:\mathbb{Q}(-1)$$ is isomorphic to $(-1)$-twist of
$\mathbb{Q}(0)$. By a $(-1)$-twist, we mean that for all $(X,P) \in
\mathcal{M}_k $ one has
$$Hom_{\mathcal{M}_k}((X,P), \mathbb{Q}(-1))=id_{Spec(k)}\circ
Hom^{-1}(X, Spec(k))\circ P \text{.}$$ This motive $\mathbb{Q}(-1)$
is called the Lefschetz motive. Now let
$$\mathbb{Q}(a):=\mathbb{Q}(-1)^{\otimes-a}$$ for $a<0$, and let $X\in
\mathcal{V}_k$ be connected. Then
$$Hom_{\mathcal{M}_k}((X,P)\otimes\mathbb{Q}(a),(Y,Q)\otimes\mathbb{Q}(b))=Q\circ CH^{dimX-a+b}(X\times Y)\circ P\text{.}$$
Similar to $h(\mathbb{P}_k^1)$, one can give a decomposition of
$h(X)$ for any $X \in \mathcal{V}_k$. Let $x\in X$ be a closed point
with residue field $\kappa (x)$. Consider the degree one cycles
$$\alpha:= deg(\kappa (x) / k)^{-1}[\{x\} \times Spec(k)] \in Hom_{\mathcal{M}_k}(h(X),
\mathbb{Q}(0))$$ and $$\beta:= [Spec(k)\times X]\in
Hom_{\mathcal{M}_k}(\mathbb{Q}(0),h(X))\text{.}$$ Then we get $
\alpha \circ \beta = id_{\mathbb{Q}(0)}$. Therefore
$id_X=id_{\mathbb{Q}(0)} + (id_X - id_{\mathbb{Q}(0)})$ is an
orthogonal decomposition where both summands are projectors. This
yields a decomposition
\begin{equation}\label{decomp}h(X)=\mathbb{Q}(0) \oplus
\tilde{h}(X)\end{equation} in $\mathcal{M}_k$, where $\tilde{h}(X)$
denotes the motive $(X,id_X - id_{\mathbb{Q}(0)})$. Note that this
decomposition in general depends on $x$. However the decomposition
$h(\mathbb{P}_k^1)= \mathbb{Q}(0) \oplus \mathbb{Q}(-1)$ is
canonical since the decomposition \ref{P1decomp} of
$\Delta_{\mathbb{P}_k^1}$ in $CH^1(\mathbb{P}_k^1 \times
\mathbb{P}_k^1)$ is independent of the choice of $x\in
\mathbb{P}_k^1(k)$.

     Let $K_0(\mathcal{M}_k)$ denote the
Grothendieck ring of the category $\mathcal{M}_k$, this is the group
generated by the isomorphism classes of motives with the relation
$$[M\oplus N]-[M] - [N] $$ for all $M,N \in \mathcal{M}_k$. The product is given by
$$[M]\cdot[N]:= [M\otimes N] \text{.}$$

    Let $X \in \mathcal{V}_k$, and let $Z$ be a closed smooth subvariety
of $X$. Then one has the following canonical isomorphism of motives
$$h(Bl_Z X) \oplus h(Z) \cong h(X)\oplus h(E)$$ where $Bl_Z X$ is
the blow up of $X$ at $Z$, and $E$ is the exceptional divisor
\cite[\S 9]{Manin2}. Hence
$$[h(Bl_Z X)]-[h(E)]=[h(X)]-[h(Z)] $$ in $K_0(\mathcal{M}_k)$. By Remark \ref{Bittner_1},
we know also that $K_0(Var_k)$ is generated by smooth projective
varieties with the blow up relations. This implies that the functor
$h$ induces a ring homomorphism $\chi \text{ : } K_0(Var_k)
\rightarrow K_0(\mathcal{M}_k)$ given by $\chi([X])=[h(X)]$ for $X
\in \mathcal{V}_k$. Note that one has
$$\chi(\mathbb{L})=[h(\mathbb{P}^1)] - [\mathbb{Q}(0)]=[\mathbb{Q}(-1)] \text{.}$$

    Let us now give the definition for the Chow groups of a motive $M\in
\mathcal{M}_k$ by
$$CH^i(M):=Hom_{\mathcal{M}_k}(\mathbb{Q}(-i),M)\text{, and }
CH_i(M):=Hom_{\mathcal{M}_k}(M,\mathbb{Q}(-i))$$ for $i\geq0$ and
$CH^i(M)=0=CH_i(M)$ for $i<0$. Note that if $M=h(X)$ for a $X\in
\mathcal{V}_k$ then we get simply $CH^i(M)=CH^i(X)$ and
$CH_i(M)=CH_i(X)$. \\ Consider a field extension $k\subset L$. There
is a functor between $\mathcal{M}_k$ and $\mathcal{M}_L$ given by
base change:
\begin{eqnarray*}
\times_k L \text{ : }\mathcal{M}_k &\longrightarrow & \mathcal{M}_L \\
(X,P) &\longmapsto& (X\times_k L, P\times_k L) \text{.}
\end{eqnarray*} In his article \cite{Chatz}, using the argumentation of \cite{BS}, Chatzistamatiou proves the following
proposition for the motives of which the degree zero Chow groups are
zero:
\begin{proposition} \label{Andre.1.}{\rm\cite[Proposition 1.2]{Chatz}}
Let $k$ be a perfect field, and $X\in \mathcal{V}_k$ be connected.
\begin{enumerate}
  \item A motive $M=(X,P)$ can be written as $M\cong N\otimes\mathbb{Q}(-1)$
  with some motive $N$ if and only if $CH_0(M\times_k L)=0$ for some field
  extension $L$ of the function field $k(X)$ of $X$.
  \item There exists an isomorphism $M\cong N \otimes \mathbb{Q}(a)$ with some
  motive $N$ and $a<0$ if and only if $CH_i(M\times_k L)=0$ for all $i<-a$ and
  all field extensions $k\subset L$.
\end{enumerate}
\end{proposition}
     We now explain that the following theorem, which is due to A.
Roitman \cite{Roitm}, implies for $X$ a hypersurface of degree
$d\leq n$ that one gets $CH_0(\tilde{h}(X))=0$.
\begin{thm}\label{Roitman}{\rm{\cite[Theorem 2]{Roitm}}}
Let $k$ be an algebraically closed field, let $X \subset
\mathbb{P}_k^{n}$ be a hypersurface of degree $d$ with $d\leq n$.
Then the subgroup of $CH_0(X)$ of degree zero cycles over
$\mathbb{Z}$, denoted by $CH^0_0(X)$ is zero.
\end{thm}
\noindent Hence the degree map
\begin{eqnarray*}
CH_0(X)\otimes \mathbb{Q} &\longrightarrow &  \mathbb{Q } \\
\sum_i n_i [Z_i] &\longmapsto & \sum_i n_i
\end{eqnarray*}
is an isomorphism.

    Note that the degree map given above is still an isomorphism over an
algebraically non-closed field (cf. \cite[Lecture 1, Appendix, Lemma
3]{Bloch}). Let $k$ be any field. Recall that
$$CH_0(X_{\overline{k}})\otimes
\mathbb{Q}=\displaystyle{\lim_{\rightarrow}} \hspace{.1cm}
CH_0(X_E)\otimes \mathbb{Q} \text{,}$$ where the direct limit is
taken over all finite extensions of $k$, and  where $\overline{k}$
denotes the algebraic closure of $k$. Let $\pi : X_E \rightarrow X$
be the base change map of a hypersurface $X$ as given in Theorem
\ref{Roitman}, over a finite extension $E $ of $k$. Then one has the
pull-back
 $\pi^* : CH_0(X) \rightarrow CH_0(X_E)$ and the push-forward
homomorphism $\pi_* : CH_0(X_E)\rightarrow CH_0(X)$, and one gets
that $\pi_* \circ \pi^*$ is multiplication by $[E:k]$ since $E / k$
is finite. Therefore after tensoring with $\mathbb{Q}$, $\pi_* \circ
\pi^*$, hence $\pi^*$ becomes injective. By passing to the direct
limits, this gives an injection $CH_0(X)\otimes \mathbb{Q}
\hookrightarrow CH_0(X_{\overline{k}})\otimes \mathbb{Q}$. Therefore
one gets that the degree map $CH_0(X)\otimes \mathbb{Q}\rightarrow
\mathbb{Q}$ is an isomorphism, since the degree map over the
algebraic closure is an isomorphism by Theorem \ref{Roitman}.

    Now let $X\subset \mathbb{P}_k^n$ be a hypersurface of degree $d
\leq n$ over a perfect field $k$, not necessarily algebraically
closed. Then one has
\begin{eqnarray*}
\mathbb{Q} \cong CH_0(X) & = & Hom_{\mathcal{M}_k}(h(X), \mathbb{Q}(0)) \\
& = &  Hom_{\mathcal{M}_k}(\mathbb{Q}(0) \oplus \tilde{h}(X), \mathbb{Q}(0)) \\
& = & \mathbb{Q}\oplus CH_0(\tilde{h}(X)) \text{,}
\end{eqnarray*}
hence $CH_0(\tilde{h}(X))=0$. Therefore applying Proposition
\ref{Andre.1.} to the motive $\tilde{h}(X)$ of Equation
\ref{decomp}, one gets the following decomposition \begin{equation}
\label{motivic_decomp} h(X)=\mathbb{Q}(0)\oplus N \otimes
\mathbb{Q}(-1)\end{equation} with some motive $N\in \mathcal{M}_k$.
\begin{remark}\label{integer_coeff}
{\rm If one considers the category of motives over $\mathbb{Z}$,
i.e. if the correspondences are Chow groups with integer
coefficients, then it is not known if having $CH_0(\tilde{h}(X))=0$
implies a decomposition of the form
\begin{equation}
\label{integral_motivic_decomp} h(X)=\mathbb{Z}(0)\oplus N \otimes
\mathbb{Z}(-1)\end{equation} with some integral motive $N$.}
 \end{remark}
It is also unknown if one in general gets a corresponding
decomposition to \ref{integral_motivic_decomp}, of the class of $X$
in the Grothendieck ring of varieties, note that this is again an
integral question. Our main concern will be to search for such
decomposition of  the classes of hypersurfaces of degree $d\leq n$
in $K_0(Var_k)$. Let us formally ask our main question which is due
to H. Esnault. Let $[X]$ denote the equivalence class of $X$ in
$K_0(Var_k)$ from now on, unless otherwise stated.

\begin{question} \label{question}
For a projective hypersurface $X \subset \mathbb{P}^{n}_k$ of degree
$d \leq n$, does one have $X(k) \neq \emptyset$ if and only if
$[X]\equiv 1 \textnormal{ mod } \mathbb{L}$ in $K_0(Var_k)$?
\end{question}
    Now clearly the spirit of Question \ref{question} varies with the
degree of $X$ and also with the base field $k$. Recall that for a
field $k$ to be $C1$ means that any homogeneous polynomial
$F(x_0,\dots,x_n)\in k[x_0,\dots,x_n]$ of degree $d\leq n$ has a
nontrivial solution over $k$. Thus over a $C1$ field Question
\ref{question} takes the following form
\begin{question} Let $k$ be a $C1$ field, and
let $X$ be as in Question \ref{question}.
Is $[X]\equiv 1 \textnormal{ mod } \mathbb{L}$?
\end{question}
In other words, having a positive answer to Question \ref{question}
over a $C1$ field $k$ implies that for \emph{every} projective
hypersurface $X \subset \mathbb{P}_k^n$ of degree $d\leq n$, one has
$[X] \equiv 1 \text{ mod } \mathbb{L}$ in $K_0(Var_k)$.

     In the remaining sections we study some particular cases which
give a positive answer to Question \ref{question}, using elementary
geometric methods. We show that for the union of $d$ hyperplanes in
$\mathbb{P}^{n}_k$ over any field $k$ with $d \leq n$, and for
quadrics over a field of characteristic zero the answer is
affirmative. In the case of a nonsingular cubic, the setting is
already far more complicated: recall that a variety $X$ of dimension
$m$ over a field $k$ is called unirational when there is a rational
map $\varphi: \mathbb{P}_k^{m}\dashrightarrow X$ such that
$\varphi(\mathbb{P}_k^{m})$ is dense in $X$ and the function field
$k(\mathbb{P}_k^{m})$ is a separable extension of $k(X)$. It is well
known that all smooth cubic hypersurfaces over an algebraically
closed field are unirational. More generally, for any cubic
hypersurface which is not a cone over a smaller dimensional cubic,
J. Koll\'{a}r \cite[Theorem 1.2]{Ko2} proved that having a
$k$-rational point is equivalent to being unirational over $k$. The
nonsingular cubic hypersurfaces in $\mathbb{P}^4_{\mathbb{C}}$,
which are in particular unirational, are proven by C. H. Clemens and
P. Griffiths to be nonrational varieties \cite[Theorem 13.12]{C-G}.
However, it is not known whether it is stably birational or not,
therefore it would be interesting to study this particular instance.
Unfortunately we are not able to answer Question \ref{question} for
this difficult case. In the higher dimensions questions of the
rationality and stable rationality of nonsingular cubic
hypersurfaces are open in general.

\newpage
\section {The Class of a Quadric in $K_0(Var_k)$} \label{ikinci_derece}  In this section we will
consider the quadrics, i.e. degree two hypersurfaces, since these
hypersurfaces would be the next natural example after hyperplanes to
look for the answer to Question \ref{question}.
\begin{thm}\label{Kollar_quadric} {\rm{(\cite[Theorem
1.11]{KoSm})}}
Let $k$ be a field of characteristic zero. Let $X \subset
\mathbb{P}_k^n$ be an irreducible quadric hypersurface. Then $X$ is
rational if and only if it has a smooth $k$-rational point.
\end{thm}
By Remark \ref{remark_ll2}, for a smooth quadric $X$, we have $X$
stably rational if and only if $[X]\equiv 1$ modulo $\mathbb{L}$
over a field of characteristic zero. This together with Theorem
\ref{Kollar_quadric} and Proposition \ref{proposition_ll3} give the
positive answer to Question \ref{question}: for a smooth quadric
$X$, one has $X(k)\neq \emptyset$ if and only if $[X] \equiv 1$
modulo $\mathbb{L}$. We like to give a proof to show that for a not
necessarily smooth quadric hypersurface over a characteristic zero
field, the answer to the Question \ref{question} is positive. We
will use the following well known lemma from linear algebra, which
is about the nondegenerate quadratic forms with a nontrivial
solution over the base field, first for writing down such an affine
fibration for a smooth quadric, and then later in Section
\ref{section_2_quadrics} for proving Theorem \ref{2quadrics}.
\begin{lemma} \label{lemma1}
Let $k$ be a field with $\text{char}(k) = 0$, and let $X=V(q)$ be a
smooth quadric hypersurface in $\mathbb{P}^n_k$. Let $x\in X(k)$ be
a point. Then we can choose coordinates in $\mathbb{P}_k^n$ such
that $$x_i(x)=0 \text{, for all } i\neq 1 \text{,}$$ and
$$q(x_0,...,x_n)= x_0x_1+q'(x_2,\dots,x_n)\text{,}$$ where $q' \in
k[x_2,\dots,x_n]$ is a quadratic form.
\end{lemma}
\begin{remark}{ \rm
With the new coordinates from Lemma \ref{lemma1}, the point $x$
corresponds to the point $[0:1:0:\dots:0]$. This will be a
convenient choice for Proposition \ref{2quadrics}.}
\end{remark}

\begin{thm} \label{quadric}
Let $k$ be a field with $\text{char}(k) = 0$.
Let $X \subset \mathbb{P}_k^n$ be a hypersurface of degree 2, then
 $[X]\equiv 1 \textnormal{ mod } \mathbb{L}$ in $K_0(Var_k)$ if and only if $X(k) \neq \emptyset$.
\end{thm}

\begin{proof}
In the case that $X$ is singular, it always has a $k$-rational
point. If $[X]\equiv 1 \textnormal{ mod } \mathbb{L}$ for $X$
smooth, then we know by Proposition \ref{proposition_ll3} that
$X(k)\neq \emptyset$.

    Now let us consider a quadric $X=V(q)$ where $q$ is a quadratic
form. First we consider the case $X$ is smooth, that is $q$ is
nondegenerate. We make the change of coordinates as given in Lemma
\ref{lemma1}, and denote by $x_0, \dots, x_n$ these new coordinates.
\noindent Let $U_0:=\{ [x_0: \dots:x_n] \in \mathbb{P}_k^n \text{ :
} x_0 \neq 0\} \subset \mathbb{P}^n_k$, then
$$U_0 \cap X = V(\frac{x_1}{x_0}+q'(\frac{x_2}{x_0},
\dots,\frac{x_n}{x_0}))  \cong \mathbb{A}^{n-1}_k\text{.}$$ Let
$Y:=V(q') \subset \mathbb{P}_k^{n-2}$. Now let us consider the
Zariski locally trivial affine fibration
\begin{eqnarray*}
p\text{ : }\mathbb{P}_k^n \setminus (U_0 \cup \{[0:1:0:\dots:0]\}) &\longrightarrow& \hspace{.9cm}\mathbb{P}_k^{n-2} \\
{[x_0:\dots:x_n]}       \hspace{0.8cm}                                            & \longmapsto & [x_2:\dots:x_n] \text{.}
\end{eqnarray*}
Then the following is a cartesian diagram:
$$
\xymatrix{ \mathbb{P}_k^n \setminus (U_0 \cup \{[0:1:0:\dots:0]\})
\ar[r] &\mathbb{P}_k^{n-2} \\
X\setminus ((X \cap U_0)\cup \{[0:1:0:\dots:0]\}) \ar[r]
\ar@{^{(}->}[u] & Y \ar@{^{(}->}[u] }
$$
Hence
$$p|_X \text{ : }X \setminus ((X \cap U_0)\cup \{[0:1:0:\dots:0]\}) \longrightarrow Y$$
is also Zariski locally trivial with fibres isomorphic to
$\mathbb{A}_k^1$. Then by Remark \ref{trivial_fibration} we have
\begin{eqnarray*}
[X] &=& [X \setminus ((X \cap U_0)\cup \{[0:1:0:\dots:0]\})] + \\ &&[(X \cap U_0)\cup \{[0:1:0:\dots:0]\}] \\
&=& \mathbb{L}\cdot[Y]+[X \cap U_0]+1 \\
&=& 1+\mathbb{L}\cdot[Y]+\mathbb{L}^{n-1}
\end{eqnarray*}
\noindent and $[X] \equiv 1\text{ mod } \mathbb{L}$ in $K_0(Var_k)$.

     Now let us consider the case where $X$ is singular. Since one can
diagonalize every quadratic form, after a change of coordinates we
can write $q(x_0,...,x_n)=a_0x_0^2+ \dots +a_rx_r^2$, with $r<n$
since $q$ is degenerate. Thus we have
$$\mathbf{P}:=\{[0:\dots:0:x_{r+1}:\dots:x_{n}] \text{ }| \text{ } [x_{r+1}:\dots: x_n] \in \mathbb{P}_k^{n-r-1}  \}  \subset X=V(q)$$
Let $Y:=V(a_0x_0^2+ \dots+ a_rx_r^2) \subset \mathbb{P}_k^r$.
Consider the projection
\begin{eqnarray*}
p: \mathbb{P}_k^n\setminus \mathbf{P} &\longrightarrow&
\mathbb{P}_k^r \\
{[x_0:\dots:x_n]} &\longmapsto&{[x_0: \dots : x_r]} \text{.}
\end{eqnarray*}
Observe that $x\in X\setminus \mathbf{P}$ if and only if $p(x)\in
Y$: indeed $p(x)\in Y$ if and only if $x \not \in \mathbf{P}$ and
$q(x)=0$. Therefore the following is a cartesian diagram
$$
\xymatrix{\mathbb{P}_k^n\setminus \mathbf{P} \ar[r]
&\mathbb{P}_k^r  \\
X \setminus \mathbf{P} \ar[r] \ar@{^{(}->}[u] & Y \ar@{^{(}->}[u] }
$$
and the projection $p$ is Zariski locally trivial, hence we see that
the map $X\setminus \mathbf{P}$ is also Zariski locally trivial
fibration with fibres isomorphic to $\mathbb{A}_k^{n-r}$, and by
Remark \ref{trivial_fibration} we obtain
\begin{eqnarray*}
[X] & = & [X \setminus \mathbf{P} ] + [\mathbf{P}] \\
& = & \mathbb{L}^{n-r}\cdot [Y] + [\mathbb{P}_k^{n-r-1}] \\
& = & 1+ \mathbb{L}+ \dots +
{\mathbb{L}}^{n-r-1}+{\mathbb{L}}^{n-r}\cdot[Y] \text{.}
\end{eqnarray*}
\noindent
Thus $[X]\equiv 1 \text{ mod } \mathbb{L}$ in $K_0(Var_k)$.
\end{proof}
\begin{corollary}
For quadric hypersurfaces Question \ref{question} has a positive
answer.
\end{corollary}
\section {The Class of a Union of Hyperplanes} \label{union_hyp}
 For a hypersurface of degree greater than two, the simplest
case to consider is the union of hyperplanes. Such a variety always
has $k$-rational points. Therefore in this case Question
\ref{question} asks if the class $[X]$ is always equivalent to $1$
modulo $\mathbb{L}$.
\begin{thm} \label{hyp1}
Let $k$ be a field, and let $X:=V(h_1 \cdots h_d) \subset
\mathbb{P}_k^n$ be a hypersurface of degree $d\leq n$ where $h_i \in
k[x_0,\dots,x_n]$ are homogeneous polynomials of degree $1$, for
$1\leq i \leq d $, then $[X]\equiv 1 \text{ mod } \mathbb{L}$ in
$K_0(Var_k)$.
\end{thm}
\begin{proof}
Let $Y:=\bigcap_{i=1}^dV(h_i)$, let $r$ be the codimension of $Y$ in
$\mathbb{P}_k^n$. Since $Y$ is an intersection of hyperplanes, it is
isomorphic to the projective space $\mathbb{P}_k^{n-r}$. Note that
$n-r \geq 0$ since we have $r \leq d \leq n$. Now we claim the
following: Let $f:=h_1\cdots h_d$, then $f\in k[x_0,
\dots,x_{r-1}]$. Indeed, after a coordinate change we get $$Y=
\bigcap_{i=1}^dV(h_i)=\{ [0:\dots:0:x_r: \dots :x_n] \text{ } |
\text{ }[x_r: \dots :x_n]\in \mathbb{P}_k^{n-r} \} \cong
\mathbb{P}_k^{n-r} \text{.}$$ Since $h_i(x)=0$ for all $x \in Y$ and
for all $i=1,\dots,d$, we find that $h_i \in k[x_0, \dots,x_{r-1}] $
for all $i=1,\dots,d$. Thus we get $f \in k[x_0, \dots,x_{r-1}]$.
Let $Z:= V(f) \subset \mathbb{P}_k^{r-1}$, and let us consider the
Zariski locally trivial affine fibration
\begin{eqnarray*}
p \text{ : } \hspace{0.3cm}\mathbb{P}_k^n \setminus Y \hspace{0.6cm}&\longrightarrow& \hspace{.9cm}\mathbb{P}_k^{r-1} \\
{[x_0: \dots:x_n]} &\longmapsto & [x_0:\dots:x_{r-1}] \text{.}
\end{eqnarray*}
Note that then we have $x\in X\setminus Y$ if and only if $p(x)\in
Z$. Hence the following is a cartesian diagram:
$$
\xymatrix{ \mathbb{P}_k^n \setminus Y  \ar[r]^{p} &
\mathbb{P}_k^{r-1} \\
X \setminus Y \ar[r] \ar@{^{(}->}[u] & Z \ar@{^{(}->}[u]}
$$
Thus
$$p|_X \text{ : }X \setminus Y \longrightarrow Z$$
is also Zariski locally trivial with fibres isomorphic to $\mathbb{A}_k^{n-r+1}$. Thus
we get $$[X \setminus Y]=\mathbb{L}^{n-r+1}\cdot[Z] \text{.}$$ Hence
by Remark \ref{trivial_fibration} we obtain
\begin{eqnarray}
[X] & = & [X\setminus Y]+[Y] \nonumber  \\
& = & \mathbb{L}^{n-r+1} \cdot [Z]+ [\mathbb{P}_k^{n-r}] \nonumber  \\
& = & 1+ \mathbb{L}+ \dots +\mathbb{L}^{n-r}+\mathbb{L}^{n-r+1}\cdot[Z] \nonumber
\end{eqnarray}
and $[X]\equiv 1 \text{ mod } \mathbb{L}$ in $K_0(Var_k)$.
\end{proof}
\begin{corollary}
Let $X$ be as in Theorem \ref{hyp1}, then $[X]\equiv 1 \text{ mod }
\mathbb{L}$ in $K_0(Var_k)$ if and only if $X(k)\neq \emptyset$.
\end{corollary}
\begin{proof}
This follows from the fact that any point of the form
$$[0:\dots:0:x_r:\dots:x_n]\in \mathbb{P}_k^n$$ is a point of $X$.
Hence $X$ has always $k$-rational points.
\end{proof}
    After checking union of hyperplanes, next example we
consider will be the hypersurfaces which become a union of
hyperplanes after a finite Galois base change.
\begin{thm} \label{hyp2}
Let $L$ be a finite Galois extension of $k$, let $X$ be a hypersurface
in $\mathbb{P}_k^n$ of degree $d \leq n$ such that $X\times_k L$ is a
union of $d$ hyperplanes over $L$, that is $X\times_k L=V(h_1\cdots h_d)$
where $h_i \in L[x_0,\dots,x_n]$ are homogeneous of degree $1$ for
$1\leq i\leq d$. Then $[X]\equiv 1 \text{ mod } \mathbb{L}$ in $K_0(Var_k)$, and $X(k) \neq \emptyset$.
\end{thm}
\begin{proof}
Let $Y:=\bigcap_{i=1}^dV(h_i)$ in $\mathbb{P}_L^n$, and let
$G:=Gal(L/k)$. Then $G$ acts on $X\times_k L$, and also on $Y$ since
one has $g\cdot y \in Y$ for all $y\in Y$. Therefore we see that
$Y/G \subset X$ is a closed immersion.

    \emph{Claim. }$Y/G$ is a linear subspace of $\mathbb{P}_k^n$. \\
In order to prove this claim, let us consider the $L$-vector
subspace $$V:=L<h_1, \dots ,h_d> \subset L^{n+1}\text{.}$$ Since $G$
acts on $V(h_1\cdots h_d)=X\times_k L$, it also acts on $V$.
Therefore one can use Lemma \ref{lemmaH90} below to deduce that
there exists homogeneous polynomials $\{h'_i\}_{i=1}^r$ in
$k[x_0,\dots,x_n]$ of degree $1$ such that $L<h'_1, \dots, h'_r>=V$,
where $r:=dim_L V$. Hence, $Y/G=\bigcap_{i=1}^r V(h'_i)$ and
$Y=(Y/G) \times_k L$. Now that $Y/G$ is a linear subspace of
$\mathbb{P}_k^n$ of codimension $r$, we can choose the coordinates
so that we have
$$Y/G=\{[0:\dots:0:x_r:\dots:x_n]\text{ } | \text{ } [x_r:\dots:x_n]
\in \mathbb{P}_k^{n-r} \} \text{.}$$ Consider the projection map
\begin{eqnarray*}
p \text{ : } \hspace{0.4cm} \mathbb{P}_k^{n}\setminus (Y/G) &\longrightarrow& \hspace{.9cm}\mathbb{P}_k^{r-1} \\
{[x_0: \dots:x_n]} &\longmapsto & [x_0:\dots:x_{r-1}] \text{.}
\end{eqnarray*}
Observe that  $$h_i \in L[x_0,\dots,x_{r-1}]$$ and
$$f:=\prod_{i=1}^d h_i \in k[x_0,\dots,x_{r-1}]$$ in a similar
manner to the polynomials $h_i$ of Theorem \ref{hyp1} (see proof of
Theorem \ref{hyp1}). Let $Z:=V(f) \subset \mathbb{P}_k^{r-1}$. With
this setting we obtain that $x\in X \setminus (Y/G)$ if and only if
$p(x) \in Z$. Therefore the following diagram
$$
\xymatrix{ \mathbb{P}_k^n \setminus (Y/G) \ar[r]^{p}
&\mathbb{P}_k^{r-1} \\
X \setminus (Y/G) \ar[r] \ar@{^{(}->}[u] & Z \ar@{^{(}->}[u] }
$$
is cartesian, $p|_X \text{ : }X \setminus (Y/G) \longrightarrow Z$
is Zariski locally trivial affine fibration with fibres isomorphic
to $\mathbb{A}_k^{n-r+1}$. Thus we are now again in the same
situation as in Theorem \ref{hyp1}, i.e, $[X]\equiv 1 \text{ mod }
\mathbb{L}$ in $K_0(Var_k)$. Moreover $X(k) \neq \emptyset$ since $
\emptyset \neq (Y/G)(k)\subset X$ is a closed immersion.
\end{proof}
\begin{lemma}[Hilbert 90] \label{lemmaH90}
Let $k$ be a field, and let $L$ be a finite Galois extension of $k$.
We denote by $G:=Gal(L/k)$ the Galois group of the extension $L/k$.
Let $V\subset L^{n}$ be a $G$- invariant $L$-vector subspace, i.e.,
$g(V) \subset V$ for all $g\in G$, then there exists a $k$-vector space
$V' \subset k^n$ such that $V'\otimes_k L=V$ where $V'\otimes_k L$ is
considered as an $L$-vector subspace of $L^n$ by the extension of the scalars.
\end{lemma}
\begin{proof}
Let $r:={dim}_L V$, and let $\{t_i\}_{i=1}^r$  be an $L$-basis of
$V$. Since $g(V)\subset V$ for all $g\in G$, $G$ acts on $V$, as
below
$$g\cdot t_i = \sum_{j=1}^r m_{ji}(g)t_{j}$$ where $m_{ij}(g)\in L$ not
all simultaneously zero, for all $t_i$. We need to prove that there
exists an $L$-basis $\{t'_i\}_{i=1}^r$ of $V$ such that $g\cdot
t'_i=t'_i$ for all $i$, because this means that $\{t'_i\}_{i=1}^r$
consist a $k$-basis for a vector space $V'$ such that $V'\otimes_k
L=V$. For this aim, let us define the map
\begin{eqnarray*}
\alpha \text{ : } \hspace{.4cm} G \hspace{.1cm} &\longrightarrow&  \hspace{.5cm} GL_r(L)\\
g \hspace{.2cm}&\longmapsto & \alpha_g := (m_{ji}(g))^{-1}_{i,j} \text{ .}
\end{eqnarray*}
Recall that a $1$-cocycle with values in $GL(V)$ is a map
$\rho \text{ : } G \longrightarrow GL(V)$ such that $\rho(g_1g_2)=\rho(g_1)g_1(\rho(g_2))$ for all $g_1 \text{, }g_2 \in G$.\\ \\
Claim: $\alpha \text{ : } G \longrightarrow GL_r(L)$ is a $1$-cocycle.
Indeed, for all $g_1 \text{, }g_2 \in G$ we have
$$g_1g_2 \cdot t_i= g_1\cdot (g_2 \cdot t_i) =\sum_{j=1}^r \sum_{k=1}^r g_1(m_{ki}(g_2))m_{jk}(g_1)t_j$$
\noindent Thus we have
\begin{eqnarray*}
\alpha_{g_1g_2} &=& ((m_{ji}(g_1g_2))_{i,j})^{-1}\\
&=&[(g_1(m_{ji}(g_2))(m_{ji}(g_1)))_{i,j}]^{-1} \\
&=& ((m_{ji}(g_1))_{i,j})^{-1}(g_1(m_{ji}(g_2)_{i,j})^{-1} \\
&=& \alpha_{g_1} g_1(\alpha_{g_2}) \text{  .}
\end{eqnarray*} \\
\noindent By the Hilbert $90$ theorem \cite[Chapter $X$, Proposition
$3$]{Serre}, we know that $$H^1(G,GL_r(L))=\{1\}\text{,}$$ i.e.,
every $1$-cocycle with values in $GL_r(L)$ is cohomologous to the
trivial $1$-cocycle that maps every element of $G$ to $id \in
GL_r(L)$. Hence $\alpha$ is cohomologous to $1$, i.e. there exists a
$B \in GL_r(L)$ such that $\alpha_g = Bg(B^{-1})$ for all $g \in G$.
Now the following claim will complete the proof of the lemma: \\
\noindent
Claim: $\{B^{-1}\cdot t_i\}_{i=1}^r$ is an $L$-basis of $V$ such that $g(B^{-1}\cdot t_i)=B^{-1}\cdot t_i$ for all $i=1,\dots ,r$.\\
Indeed, we have \begin{eqnarray*}
g(B^{-1}\cdot t_i) &=& g(B^{-1})g \cdot t_i \\
&=& g(B)^{-1} \alpha_g^{-1} \cdot t_i \\
&=& g(B)^{-1} g(B)B^{-1}\cdot t_i \\
&=& B^{-1}\cdot t_i \text{   .}
        \end{eqnarray*}
\end{proof}
\begin{corollary}
For hypersurfaces of the type given in Theorem \ref{hyp1} and
\ref{hyp2}, Question \ref{question} has a positive answer.
\end{corollary}
\section {Cubic Hypersurfaces} \label{kubik}  In this section we will
consider the case of degree three hypersurfaces for Question
\ref{question}. Therefore we are only interested in the
hypersurfaces that live in at least three dimensional projective
space. Note that over an algebraically closed field $k$, it is known
that a smooth cubic surface is rational \cite[Chapter IV, Theorem
24.1]{Manin1}. For a smooth projective hypersurface $X$ over a field
of characteristic zero, (stable) rationality implies that the class
$[X]$ is equivalent to $1$ modulo $\mathbb{L}$ in $K_0(Var_k)$, by
Proposition \ref{larsen_lunts_2}. On the other hand there are smooth
cubic hypersurfaces that are not rational, like smooth cubic
threefolds which are proven to be non-rational by Clemens and
Griffiths in their paper \cite{C-G}. It is still unknown whether
smooth cubic threefolds are stably rational or not. If an
irreducible projective cubic hypersurface that is not a cone over a
cubic of lower dimension has a $k$-rational singular point, then it
is rational \cite[Chapter 1, Section 5, Example 1.28]{KoSm}. However
in the singular case rational varieties do not necessarily have
class $1$ modulo $\mathbb{L}$, but it is possible to show that a
singular cubic hypersurface with a rational singular point, actually
has this property. This is what we are going to prove next.
\begin{thm} \label{cubicsingular}
Let $k$ be a $C1$ field, and let $X$ be a hypersurface of degree $3$
in $\mathbb{P}_k^n$. Denote by $X_{sing}(k)$ the set of $k$-rational
points of the singular locus of $X$. If $X_{sing}(k) \neq
\emptyset$, then $[X]\equiv 1 \text{ mod } \mathbb{L}$ in
$K_0(Var_k)$.
\end{thm}
\begin{proof}
Let $x \in X_{sing}(k)$. After a change of coordinates, we get
$x=[0:\dots:0:1]$. Let $X=V(f)$. Since $x=[0:\dots:0:1]\in X(k)$,
$f$ has the following form with these new coordinates
$$f(x_0,\dots, x_n)=x_n^2f_1(x_0,\dots,x_{n-1})+x_nf_2(x_0,\dots,x_{n-1})+f_3(x_0,\dots,x_{n-1})$$
where $f_i$ are homogeneous polynomials of degree $i$, $i=1,2,3$.
Moreover $x$ is a singular point of $X$, i.e.
$${\frac {\partial f}{\partial x_i}}|_{x}=0, \text{ for all } 0\leq i\leq n \text{.}  $$
Let us note that we have in particular $${\frac {\partial f_1}{\partial x_i}}|_{x}=0, \text{ for all } 0\leq i\leq n-1 \text{.} $$
Hence we get $f_1=0$. Now let us consider the following Zariski
locally trivial affine fibration:
\begin{eqnarray*}
\pi \text{ : } \hspace{.5cm}\mathbb{P}^n_k \setminus \{x\} \hspace{.2cm} &\longrightarrow& \hspace{.9cm}\mathbb{P}_k^{n-1} \\
{[x_0:\dots:x_n]} &\longmapsto& [x_0:\dots:x_{n-1}] \text{.}
\end{eqnarray*}
Note that
$$\pi^{-1}(p)=\{[ p_0:\dots: p_{n-1}:\gamma] \text{ } | \text{ }\gamma \in k\} \cong \mathbb{A}_k^1 $$
and that $\pi^{-1}(p)\cup \{x\} \cong \mathbb{P}_k^1$ for all
$p=[p_0:\dots:p_{n-1}] \in \mathbb{P}_k^{n-1}$. Let us denote by
$\pi_{X}$ the restriction of $\pi$ to $X$. Therefore we get that
\begin{eqnarray} [X] &=& {[X\setminus \{x\}]}+ [\{x\}] \nonumber \\
                      &=& 1+ [\pi_X^{-1}(\mathbb{P}_k^{n-1})]
                      \label{gulus}
\end{eqnarray}
by Remark \ref{trivial_fibration}. For any $p\in \mathbb{P}_k^{n-1}$
we have  one of the following possibilities for the fibre of $\pi_X$
at $p$
$$
\pi_{X}^{-1}(p) = \left\{
                                                                  \begin{array}{l}
                                                                    y\in X\setminus \{x\} \\
                                                                    \pi^{-1}(p) \\
                                                                    \emptyset
                                                                  \end{array}
                                                                  \right.
$$
First of all, let us observe that
\begin{eqnarray*}
Y &:=& \{p \in \mathbb{P}_k^{n-1} \text{ } | \text{ }
\pi_X^{-1}(p)=y
\in X\setminus\{x\}\}\hphantom{..} \cong \hphantom{..}\pi_X^{-1}(Y) \vspace{.5cm}\\
  &=&\{p\in\mathbb{P}_k^{n-1} \text{ } | \text{ }
  \pi_X^{-1}(p)=[p_0:\dots:p_{n-1}:-f_3(p)/f_2(p)]\} \\
  &=&\{p\in\mathbb{P}_k^{n-1} \text{ } | \text{ } f_2(p)\neq 0\} \\
  &=&\mathbb{P}_k^{n-1} \setminus V(f_2) \text{.} \vspace{.5cm}
\end{eqnarray*}
Note that $k$ being a $C1$ field assures that
$$V(f_2)(k)\neq \emptyset \text{,}$$ which implies by the Theorem
\ref{quadric} that $[V(f_2)]\equiv 1 \text{ mod } \mathbb{L}$ in
$K_0(Var_k)$. Hence
\begin{eqnarray}
[\pi_X^{-1}(Y)]\hphantom{..}=\hphantom{..}[Y]&=&[\mathbb{P}_k^{n-1} \setminus V(f_2)] \nonumber \\
   &=&[\mathbb{P}_k^{n-1}] - [V(f_2)] \nonumber \\
   &\equiv & 0 \hphantom{\mathbb{P}_k^{n-1}} \text{mod } \mathbb{L}
   \label{papatya}
\end{eqnarray}

Let us denote by $Z$ the set of points of which the pre-images under
$\pi$ is completely contained in $X$:
\begin{eqnarray*}
Z &:=& \{p \in \mathbb{P}_k^{n-1} \text{ } | \text{ } \pi^{-1}(p)
\subset X\setminus\{x\}\} \\
  &=&\{p\in\mathbb{P}_k^{n-1} \text{ } | \text{ } \gamma f_2(p)+f_3(p)=0 \text{, } \forall \gamma \in k \} \\
  &=& V(f_2,f_3)
\end{eqnarray*}
Hence the diagram
$$
\xymatrix{\mathbb{P}_k^n \setminus \{x\} \ar[r]^{\pi}
&\mathbb{P}_k^{n-1} \\
\pi^{-1}(Z) \ar[r]^{\pi_X} \ar@{^{(}->}[u] & Z \ar@{^{(}->}[u] }
$$
is cartesian, and $\pi_X^{-1}(Z)\rightarrow Z$ is a Zariski locally
trivial $\mathbb{A}_k^1$-fibration. This yields by Remark
\ref{trivial_fibration} that
\begin{equation}\label{lale}
[\pi_X^{-1}(Z)]=[\pi^{-1}(Z)]= \mathbb{L} \cdot [V(f_2,f_3)]
\end{equation}
Thus we obtain
\begin{eqnarray*}
[X]&=&1+ [\pi_X^{-1}(\mathbb{P}_k^{n-1})] \hspace{.5cm}\text{by Equation \ref{gulus},} \\
   &=&1+[\pi_X^{-1}(Y)\sqcup \pi_X^{-1}(Z)] \\
   &=&1+[Y]+\mathbb{L}\cdot[V(f_2,f_3)] \hspace{.5cm}\text{by Equation \ref{lale},}\\
   &\equiv& 1 \hspace{.2cm}\text{ mod } \mathbb{L} \hspace{.5cm}\text{by Equation \ref{papatya} .}
\end{eqnarray*}
in $K_0(Var_k)$.
\end{proof}
\begin{corollary}
For cubic hypersurfaces over $k$ with a singular $k$-rational point,
Question \ref{question} is positively answered.
\end{corollary}
\section {Union of Two Quadric Hypersurfaces} \label{section_2_quadrics}  In this
section we consider a particular example of quartic hypersurfaces,
namely ones in $\mathbb{P}_k^n$ for any $n \geq 4$, which consist of
the union of two quadric hypersurfaces, one of which is smooth. We
need to assume that $k$ is algebraically closed to conclude with the
proof. We will compute the class in $K_0(Var_k)$ modulo the class
$\mathbb{L}$.
\begin{thm} \label{2quadrics}
Let $k$ be an algebraically closed field of characteristic zero, let
$X$ be a union of two quadric hypersurfaces $Q_1, Q_2 \subset
\mathbb{P}_k^n$, where $Q_1$ is smooth and $n \geq 4$. Then $[X]
\equiv 1 \text{ mod } \mathbb{L}$ in $K_0(Var_k)$.
\end{thm}
\begin{proof}
If $Q_1=Q_2$, then $[Q_1\cup Q_2]=[Q_1\cup Q_1]=[Q_1]$. Hence it
reduces to the case of Theorem \ref{quadric}, and we are done. For
distinct $Q_1, Q_2$  we have $$[X]=[Q_1\cup
Q_2]=[Q_1]+[Q_2]-[Q_1\cap Q_2]$$ in $K_0(Var_k)$. Since $k$ is
algebraically closed, both $Q_1$ and $Q_2$ have $k$-rational points.
Hence, by Theorem \ref{quadric} we get
\begin{eqnarray}\label{q_i}[Q_i] &\equiv& 1\text{ mod } \mathbb{L} \text{,    } \text{ } \text{ }i=1\text{,}2 \text{,}\end{eqnarray}
 and therefore $[X]\equiv 1 \text{ mod } \mathbb{L}$ if and only if $[Q_1\cap Q_2]\equiv 1 \text{ mod } \mathbb{L}$ in $K_0(Var_k)$.
 Let $$Q_i:= V(q_i)\text{, }\text{ }i=1,2\text{.}$$ Since we are over an algebraically closed field, the intersection $Q_1\cap Q_2$ also has a
 $k$-rational point $x$. Since $Q_1$ is smooth, we can make a change of coordinates in order to get $q_1(x_0,\dots,x_n)=x_0x_1-h(x_2,\dots,x_n)$.
 To do this, we apply Lemma \ref{lemma1} to $q_1$ and $x$ . Hence we get coordinates in which $q_1$ is of the desired form,
 and the point $x$ becomes the point $[0:1:0:\dots:0]$. Now since $[0:1:0:\dots:0]\in Q_2$, the monomial ${x_1}^2$ does not appear in $q_2$.
 Thus we can write $$q_2(x_0,\dots,x_n)=x_0L_0(x_0,x_2,\dots,x_n)+x_1L_1(x_0,x_2,\dots,x_n)+R(x_2,\dots,x_n)$$ where $L_0,L_1,R$ are
 homogeneous polynomials, $L_0,L_1$ are linear and $R$ is of degree $2$. Let $U_0:=\{[x_0:\dots:x_n] \text{ } | \text{ }x_0\neq
 0\}$. We will calculate the class of the intersection of $Q_1\cap
 Q_2$ with $U_0$. Let $Q_1^{\{x_0\neq 0 \}}:=Q_1\cap U_0=V(q_1|_{U_0}) $. We have
$$q_1|_{U_0}(\frac{x_1}{x_0},\dots,\frac{x_n}{x_0})=\frac{x_1}{x_0}-h(\frac{x_2}{x_0},\dots,\frac{x_n}{x_0})$$
and $Q_1^{\{x_0 \neq 0 \}} \cong \mathbb{A}^{n-1}_k$ where the
isomorphism is given by \begin{eqnarray*}
 \varphi \text{ : } \hspace{.4cm} Q_1^{\{x_0 \neq 0\}} & \longrightarrow & \hspace{.6cm}\mathbb{A}_k^{n-1} \\
 (\frac{x_1}{x_0}, \dots, \frac{x_n}{x_0}) &\longmapsto& (\frac{x_2}{x_0}, \dots, \frac{x_n}{x_0})
 \end{eqnarray*}
Let us now consider the intersection $Q_1^{\{x_0\neq 0 \}} \cap
Q_2$. Via the isomorphism $\varphi$, it is defined by the following
polynomial:
\begin{multline*}g(\frac{x_2}{x_0},\dots,\frac{x_n}{x_0}):=L_0(1,\frac{x_2}{x_0},\dots,\frac{x_n}{x_0})+h(\frac{x_2}{x_0},\dots,\frac{x_n}{x_0})L_1(1,\frac{x_2}{x_0},\dots,\frac{x_n}{x_0})+\\ R(\frac{x_2}{x_0},\dots,\frac{x_n}{x_0})\end{multline*}
i.e., $$\varphi (Q_1^{\{x_0 \neq 0\}} \cap Q_2)=V(g)\subset
\mathbb{A}_k^{n-1}\text{.}$$ We embed $\mathbb{A}^{n-1}_k \subset
\mathbb{P}^{n-1}_k$ as $\{y_1 \neq 0 \}$ with homogeneous
coordinates $\{y_1, \dots,y_n\}$ for $\mathbb{P}_k^{n-1}$. Then the
closure $$Y:= \overline{\varphi(Q_1^{\{x_0\neq 0 \}} \cap Q_2)}
\subset \mathbb{P}_k^{n-1}$$ is defined by the homogenization of
$g$:
$$\overline{g}(y_1,\dots,y_n)= y_1^2L_0(y_1,\dots,y_n)+h(y_2,\dots,y_n)L_1(y_1,\dots,y_n)+y_1R(y_2,\dots,y_n) \text{.}$$  Thus \\
$$ \varphi (Q_1^{\{x_0\neq 0 \}} \cap Q_2)=V(g)=Y\setminus (Y\cap V(y_1)) \text{,}$$ and we get
\begin{equation}[\varphi(Q_1^{\{x_0\neq 0 \}} \cap Q_2)] = [Y]-[Y\cap V(y_1)] \text{.}\end{equation} \\
\noindent
The intersection $Y\cap V(y_1)$ is the vanishing locus of $$h(y_2,\dots,y_n)L_1(0,y_2,\dots,y_n)\text{.}$$ Hence \begin{eqnarray*}[Y\cap V(y_1)] &=&[V(h L_1)\cap V(y_1)] \\
&=&[V(h)\cap V(y_1)]+[V(L_1) \cap V(y_1)]-[V(h)\cap V(L_1) \cap V(y_1)]\end{eqnarray*} \\ in $K_0(Var_k)$. Here $V(h)\cap V(y_1)$ is a quadric hypersurface and $V(L_1)\cap V(y_1)$ is a hyperplane in $\mathbb{P}_k^{n-2}$. Since $n-2\geq 2$ by the assumption, we get $$[V(h)\cap V(y_1)]\text{, } [V(L_1)\cap V(y_1)] \equiv 1 \text{ mod } \mathbb{L}$$ in $K_0(Var_k)$. Hence \\ \begin{equation} \label{hyperinf}
[Y\cap V(y_1)] \equiv 2 - [V(h)\cap V(L_1)\cap V(y_1)] \text{ } \text{ mod} \text{ } \mathbb{L}                                                                                                                                                                                                                                                                                                          \end{equation} \\
in $K_0(Var_k)$. Now let us examine $$Q_1^{\{x_0=0\}} \cap
Q_2:=Q_1\cap Q_2 \setminus Q_1^{\{x_0\neq 0\}} \cap Q_2 \text{.} $$
It is defined by the following polynomials:
\begin{eqnarray*}
q_1(x_0,\dots,x_n)|_{x_0=0} &=& h(x_2,\dots,x_n) \\
q_2(x_0,\dots,x_n)|_{x_0=0} &=& x_1L_1(0,x_2,\dots,x_n)+R(x_2,\dots,x_n)
\end{eqnarray*}
\noindent Now we will calculate the class of the intersection of
$Q_1 \cap Q_2$ with the complement of $U_0$. Consider the projection
map
\begin{eqnarray*}
\pi \text{ : }\hspace{.5cm} \mathbb{P}^{n-1}_k \setminus \{p\} &\longrightarrow& \hspace{.8cm} \mathbb{P}^{n-2}_k \\
{[x_1: \dots :x_n]} &\longmapsto & [x_2: \dots :x_n]
\end{eqnarray*}
where $p:=[1:0:\dots:0]$. Now this projection map induces an
isomorphism \begin{eqnarray*} (Q_1^{\{x_0=0 \}} \cap Q_2) \setminus
((V(L_1)\cap V(x_0))\cup \{p\}) & \cong &(V(h)\cap V(x_0)) \setminus
\\ &&  ((V(h)\cap V(L_1)\cap V(x_0)) \end{eqnarray*} and thus
\begin{multline} \label{iso1} [Q_1^{\{x_0=0\}}\cap Q_2 \setminus
((V(L_1)\cap V(x_0))\cup \{p\})] = \\ [V(h)\cap V(x_0)] - [V(h)\cap
V(L_1)\cap V(x_0)]
\end{multline} in $K_0(Var_k)$.\\
\noindent Since $n-2 \geq 2$, we have $[V(h)\cap V(x_0)] \equiv 1 $
mod $\mathbb{L}$ in $K_0(Var_k)$. Besides, the projection map $\pi$
induces a Zariski locally trivial $\mathbb{A}_k^1$-fibration
$$Q_1^{\{x_0=0\}}\cap Q_2\cap V(L_1)\setminus \{p\} \longrightarrow
V(L_1 ) \cap V(R)\cap V(h) \cap V(x_0) \text{.} $$  Hence we have
\begin{equation} \label{fibration} [Q_1^{\{x_0=0\}}\cap Q_2\cap
V(L_1)\setminus \{p\}]= \mathbb{L}\cdot [V(L_1 ) \cap V(R)\cap V(h)
\cap V(x_0)]
\end{equation} in $K_0(Var_k)$. By Equality \ref{iso1}  and Equality \ref{fibration}, we obtain  \begin{eqnarray}
[Q_1^{\{x_0=0\}} \cap Q_2]
&=&[Q_1^{\{x_0=0\}}\cap Q_2 \setminus ((V(L_1)\cap V(x_0))\cup \{p\})]+ \nonumber \\
&&[Q_1^{\{x_0=0\}}\cap Q_2\cap V(L_1)\setminus \{p\}]+1 \nonumber \\
&=& [V(h)\cap V(x_0)]-[V(h)\cap V(L_1)\cap V(x_0)]+ \nonumber\\
&& \mathbb{L}\cdot [V(L_1 ) \cap V(R)\cap V(h) \cap V(x_0)] +1 \nonumber \\
&\equiv &2- [V(h)\cap V(L_1)\cap V(x_0)] \text{ mod } \mathbb{L} \label{x_0=0}
\end{eqnarray}
in $K_0(Var_k)$. Here let us note that $V(L_1)\cap V(y_1)=V(L_1)\cap V(x_0)$. Therefore, putting Congruence \ref{hyperinf} and Congruence \ref{x_0=0} together, we get
\begin{eqnarray}
[Q_1\cap Q_2] &=& [Q_1^{\{x_0\neq 0 \}} \cap Q_2]+ [Q_1^{\{x_0=0\}} \cap Q_2] \nonumber \\
&=& [Y]-[Y\cap V(y_1)]+ [Q_1^{\{x_0=0\}} \cap Q_2] \nonumber \\
&\equiv& [Y] \text{ mod } \mathbb{L} \label{q1capq2}
\end{eqnarray}
in $K_0(Var_k)$. Now let us examine the class of $Y$ in
$K_0(Var_k)$. We consider the following subvariety of $Y$
\begin{multline*}S:= \{[y_1:\dots:y_n]\in \mathbb{P}_k^{n-1} \text{ } |
\text{ }y_1=h(y_2,\dots,y_n)=L_1(y_1,\dots,y_n) \\
\hphantom{,y_n)}=R(y_2,\dots,y_n)=0 \}\text{. }\end{multline*}  For
each $s:=[s_1:\dots:s_n]\in S$, we have
$$\frac{\partial \overline{g}}{\partial
y_1}(s)=2s_1L_0(s)+s_1^2\frac{\partial L_0}{\partial
y_1}(s)+h(s)\frac{\partial L_1}{\partial y_1}(s)=0\text{,}$$ and for
$2\leq i \leq n$,
$$\frac{\partial \overline{g}}{\partial y_i}(s)= s_1^2\frac{\partial L_0}{\partial y_i}(s)+ \frac{\partial h}{\partial y_i}(s) L_1(s)+ h(s)\frac{\partial L_1}{\partial y_i}(s)+s_1\frac{\partial R}{\partial y_i}(s)=0 \text{.}$$
Hence $S \subset Sing(Y)$. Now for $n\geq 5$, $S\neq \emptyset$,
therefore $Sing(Y)\neq \emptyset$, which implies that $[Y] \equiv 1$
mod $\mathbb{L}$ in $K_0(Var_k)$, by Theorem \ref{cubicsingular}. In
the case that $n=4$, $Y\subset \mathbb{P}_k^3$ is in general smooth.
However, in $\mathbb{P}_k^3$ a smooth cubic surface is always
rational \cite[Chapter IV, Theorem 24.1]{Manin1}. Thus
\begin{equation} \label{Y} [Y]\equiv 1 \text{ mod } \mathbb{L}
\text{ in } K_0(Var_k)
\end{equation} for all $n\geq 4$.
Hence we have
\begin{eqnarray*}
[X]&=& [Q_1]+[Q_2]-[Q_1\cap Q_2] \\
&=& 2-[Q_1\cap Q_2] \text{, } \hspace{1.1cm}\text{ by Congruence \ref{q_i},}  \\
&\equiv& 2-[Y] \text{ mod }  \mathbb{L} \text{, } \hspace{.9cm} \text{ by Congruence \ref{q1capq2},} \\
&\equiv& 1 \text{ mod }  \mathbb{L} \text{, } \hspace{1.8cm} \text{ by Congruence \ref{Y}}
\end{eqnarray*}
in $K_0(Var_k)$.\\
\end{proof}
\begin{corollary}
Quartic hypersurfaces of the form described in Theorem
\ref{2quadrics} gives a positive answer to Question \ref{question}.
\end{corollary}


\begin{thebibliography}{99}

\bibitem{WMKA} \textbf{Abromovich, D.; Karu, K.; Matsuki, K.; Wlodarzcyk, J.}. Torification and factorization of birational maps,
\emph{J. Algebraic Geom.}, \textbf{9}, 2000, 425-449.

\bibitem{Motif} \textbf{Andr\'{e}, Y.}. Une Introduction aux Motifs (Motifs Purs, Motifs Mixtes,
Periodes), \emph{Soc. Math. de France}, 2004.

\bibitem{BCTSSD} \textbf{Beauville, A.; Colliot-Th\'{e}l\`{e}ne, J-L.;
Sansuc J-J. and Swinnerton-Dyer P.}. Vari\'{e}t\'{e}s stablement
rationnelles non rationnelles, \emph{Ann. of Math.}, \textbf{121},
1985, 283-318.

\bibitem{F.H} \textbf{Bittner, F.}. The universal Euler characteristic for
varieties of characteristic zero, \emph{Foundation Compositio
Mathematica}, 2004, 1011-1032.

\bibitem{Bloch} \textbf{Bloch, S.}. Lectures on algebraic cycles,
\emph{Durham, NC  : Mathematics Department, Duke Univ.}, 1980.

\bibitem{BS} \textbf{Bloch, S. and Srinivas, V.}. Remarks on correspondences and algebraic cycles,
\emph{Amer. J. Math}, \textbf{105, no. 5}, 1983, 1235-1253.

\bibitem{Chatz} \textbf{Chatzistamatiou, A.}. First coniveau notch of the Dwork family and its mirror,
\emph{Math. Res. Lett.}, \textbf{16}, 2009, 563-575.

\bibitem{C-G} \textbf{Clemens, C. H. and Griffiths, P. A. }. The intermediate
Jacobian of the cubic threefold, \emph{Ann. of Math }, \textbf{Vol.
95, No. 2}, 1972, 281-356.

\bibitem{KoSm} \textbf{Corti, A.; Kollar, J.; Smith, K. E.}. Rational and
nearly rational varieties, \emph{Cambridge University Press},2004.

\bibitem{GH} \textbf{Griffiths, P. and Harris, J.}. Principles of Algebraic Geometry, \emph{Wiley-Interscience}, 1978.

\bibitem{SeCo} \textbf{Grothendieck, A.; (eds.): Colmez, P., Serre, J.-P. }.  Correspondance
Grothendieck-Serre, \emph{ Soci\'{e}t\'{e} Math\'{e}matique de
France}, \textbf{ vol. 2}, 2001.

\bibitem{Hart} \textbf{Hartshorne, R.}. Algebraic Geometry, \emph{Springer-Verlag}, 1977.

\bibitem{Ko2} \textbf{Koll\'{a}r, J.}. Unirationality of cubic
hypersurfaces, \emph{Math. Res. Lett. }, \textbf{198(1)}, 2005,
27-35.

\bibitem{Ko} \textbf{Koll\'{a}r, J.}. Conics in the Grothendieck ring,
\emph{Adv. Math., Math. Res. Lett. }, \textbf{198(1)}, 2005, 27-35.

\bibitem{LaLu} \textbf{Larsen, M. and Lunts, V. A.}.
Motivic measures and stable birational geometry, \emph{Mosc. Math.
J.}, \textbf{3(1)}, 2003, 85-95.

\bibitem{kil} \textbf{Liao, X.}. Stable birational equivalence and
geometric Chevalley-Warning, arXiv:1110.2554v1 [math.AG].

\bibitem{Manin2} \textbf{Manin, Yu. I.}. Correspondences, motives and monoidal
transformations, \emph{Math. USSR Sb. }, \textbf{6 : 4}, 439-470;
\emph{Mat. Sb.}, \textbf{77 : 4}, 1968, 475-507.

\bibitem{Manin1} \textbf{Manin, Yu. I.}. Cubic Forms. Algebra, geometry, arithmetic, \emph{North-Holland}, (Translated from Russian), 1986.

\bibitem{Nic} \textbf{Nicaise, J.}. A trace formula for varieties over a discretely valued field,
\emph{J. Reine Angew. Math.},\textbf{650}, 2011, 193-238,
arxiv:0805.1323v2.

\bibitem{Poonen} \textbf{Poonen, B.}. The Grothendieck ring of varieties is not a
domain, \emph{Math. Res. Lett. }, \textbf{9(4)}, 2002, 493-497.

\bibitem{Roitm} \textbf{Roitman, A. A.}. Rational equivalence of $0$-cycles, \emph{Mat. Zametki }, \textbf{Vol. 28, No.
1}, 1980, 85-90.

\bibitem{Scholl} \textbf{Scholl, A. J.}. Classical motives,
\emph{Motives, Seattle 1991, ed. U. Jannsen, S. Kleiman, J-P. Serre.
Proc Symp. Pure Math}, \textbf{55, part 1}, 1994, 163-187.

\bibitem{Serre} \textbf{Serre, J.-P.}. Local Fields, \emph{Berlin-New York, Springer Verlag}, 1980.


\end{thebibliography}
\end{document}